\documentclass[12pt]{article}
\usepackage{amssymb}
\usepackage{mathrsfs}
\usepackage{pst-poly}     
\usepackage{cite}
\usepackage{tikz}
\usepackage{mathrsfs}
\usepackage{amssymb}
\usepackage{amsmath}
\usepackage{bm}
\usepackage{graphicx}
\usepackage{pst-poly}  
\usepackage{pst-plot}  

\newtheorem{thm}{Theorem}[section]
\newtheorem{lem}[thm]{Lemma}

\newtheorem{cor}[thm]{Corollary}

\newenvironment{pf}[1][Proof]{\noindent\textbf{#1.} }{\hfill\rule{1mm}{2mm}}

\makeatletter \@addtoreset{equation}{section} \makeatother

\begin{document}
\title{\bf Generalized Connectivity of Star Graphs\thanks {The work was supported by
NNSF of China (No.11071233, 61272008).}}
\author
{Xiang-Jun Li$^{a,b}$
\quad Jun-Ming Xu$^a$
\footnote{Corresponding author: xujm@ustc.edu.cn (J.-M. Xu)}\\
{\small $^a$School of Mathematical Sciences}\\
{\small University of Science
and Technology of China,}\\
{\small  Wentsun Wu Key Laboratory of CAS, Hefei, 230026, China}  \\
{\small $^b$School of Information and Mathematics,}\\
{\small Yangtze University, Jingzhou, Hubei, 434023, China}\\
 }
\date{}
 \maketitle

\begin{abstract}

This paper shows that, for any integers $n$ and $k$ with $0\leqslant
k \leqslant n-2$, at least $(k+1)!(n-k-1)$ vertices or edges have to
be removed from an $n$-dimensional star graph to make it
disconnected and no vertices of degree less than $k$. The result
gives an affirmative answer to the conjecture proposed by Wan and
Zhang [Applied Mathematics Letters, 22 (2009), 264-267]. \vskip6pt



\end{abstract}

\section{Introduction}
It is well known that interconnection networks play an important
role in parallel computing/communication systems. An interconnection
network can be modeled by a graph $G=(V, E)$, where $V$ is the set
of processors and $E$ is the set of communication links in the
network.

A subset $S\subset V(G)$ (resp. $F\subset E(G)$) of a connected
graph $G$ is called a {\it vertex-cut} (resp. {\it edge-cut}) if
$G-S$ (resp. $G-F$) is disconnected. The {\it connectivity}
$\kappa(G)$ (resp. {\it edge-connectivity} $\lambda(G)$ ) of $G$ is
defined as the minimum cardinality over all vertex-cuts (resp.
edge-cuts) of $G$. The connectivity $\kappa(G)$ and
edge-connectivity $\lambda(G)$ of a graph $G$ are two important
measurements for fault tolerance of the network since the larger
$\kappa(G)$ or $\lambda(G)$ is, the more reliable the network is.
Esfahanian~\cite{E89} proposed the concept of restricted
connectivity, Latifi {\it et al.}~\cite{lhm94} generalized it to
restricted $k$-connectivity which can measure fault tolerance of an
interconnection network more accurately than the classical
connectivity. The concepts stated here are slightly different from
theirs.

A subset $S\subset V(G)$ (resp. $F\subset E(G)$) of a connected
graph $G$, if any, is called a {\it $k$-vertex-cut} (resp. {\it
edge-cut}), if $G-S$ (resp. $G-F$) is disconnected and has the
minimum degree at least $k$. The {\it $k$-super connectivity} (resp.
{\it edge-connectivity}) of $G$, denoted by $\kappa_s^{(k)}(G)$
(resp. $\lambda_s^{(k)}(G)$), is defined as the minimum cardinality
over all $k$-vertex-cuts (resp. $k$-edge-cuts) of $G$. For any graph
$G$ and any integer $k$, determining $\kappa_s^{(k)}(G)$ and
$\lambda_s^{(k)}(G)$ is quite difficult, there is no known
polynomial algorithm to compute them yet. In fact, the existence of
$\kappa_s^{(k)}(G)$ and $\lambda_s^{(k)}(G)$ is an open problem so
far when $k\geqslant 1$. Only a little knowledge of results have
been known on $\kappa_s^{(k)}$ and $\lambda_s^{(k)}$ for some
special classes of graphs for any $k$.

As an attractive alternative network to the hypercube, the
$n$-dimensional star graph $S_{n}$ is proposed by Akers {\it et
al.}~\cite{ak89}. Since it has superior degree and diameter to the
hypercube as well as it is highly hierarchical and
symmetrical~\cite{dt94}, the star graph $S_{n}$ has received
considerable attention in recent years. In particular, Cheng and
Lipman~\cite{cl02}, Hu and Yang~\cite{hy97} and Rouskov {\it et
al.}~\cite{rls96}, independently, determined
$\kappa_s^{(1)}(S_n)=2n-4$ for $n\geqslant 3$. Yang~{\it et
al.}~\cite{ylm10} proved $\lambda_s^{(2)}(S_{n})=6(n-3)$ for
$n\geqslant 4$. Wan and Zhang~\cite{wz09} showed that
$\kappa_s^{(2)}(S_{n})=6(n-3)$  for $n\geqslant 4$ and conjectured
that $\kappa_s^{(k)}(S_{n})=(k+1)!(n-k-1)$ for $ k\leqslant n-2$. In
this paper, we give an affirmative answer to the conjecture and
generalize the above-mentioned results by proving that $\kappa
_s^{(k)}(S_{n})=\lambda _s^{(k)}(S_{n}) =(k+1)!(n-k-1)$ for any $k$
with $0\leqslant k\leqslant n-2$.

In Section 2, we recall the two structures of $S_{n}$ and some
lemmas to be used in our proofs. The proof of the main results is in
Section 3. We conclude our work in Section 4.

\section{Definitions and lemmas}

For a given integer $n$ with $n\geqslant 2$, let
$I_n=\{1,2,\ldots,n\}$, $I'_n=\{2,\ldots,n\}$ and $P(n)=\{
p_{1}p_{2}\ldots p_{n}:\ p_{i}\in I_n, p_{i}\neq p_{j}, 1\leqslant
i\neq j\leqslant n\}$, the set of permutations on $I_n$. Clearly,
$|P(n)|=n\,!$. For an element $p=p_{1}\ldots p_j\ldots p_{n}\in
P(n)$, the digit $p_j$ is called the symbol in the $j$-th position
(or dimension) in $p$.

The $n$-dimensional star graph, denoted by $S_{n}$, is an undirected
graph with vertex-set $P(n)$. There is an edge between any two
vertices if and only if their labels differ only in the first and
another position. In other words, two vertices $u=p_{1}p_{2}\ldots
p_{i}\ldots p_{n}$ and $v=p'_{1}p'_{2}\ldots p'_{i}\ldots p'_{n}$
are adjacent if and only if $v=p_{i}p_{2}\ldots
p_{i-1}p_1p_{i+1}\ldots p_{n}$ for some $i\in I'_n$.

Like the hypercube, the star graph is a vertex- and edge-transitive
graph with degree $(n-1)$~\cite{ak86}. The following properties of
$S_n$ are very useful for our proof.

\begin{lem}\label{lem2.3}
\textnormal{(see Akers {\it et al.}~\cite{ak89},1989)} 
$\kappa
(S_{n})=\lambda (S_n)=n-1$ for $n\geqslant 2$.
\end{lem}

\begin{figure}[h]
\begin{center}
\begin{tikzpicture}[scale=.6]

  \tikzstyle{every node}=[draw,circle,fill=white,minimum size=5pt,
                            inner sep=0pt]
 \draw [line width = 1pt,white](-15,3.6) 
         ++(0:1.4cm) 
         ++(0:2.0cm) 
        ++(-90:6.2cm) node [label=-90:
        {\scriptsize$ \textcolor[rgb]{0.00,0.00,0.00}{ Partitioned \ along\ dimension\ 4}$}]{};

  \draw [line width = 1pt,red](-15,3.6) node (2431) [label=135:{\scriptsize$\textcolor[rgb]{0.00,0.00,0.00}{243}1$}] {}
        -- ++(0:1.4cm) node (3421) [label=45: {\scriptsize$\textcolor[rgb]{0.00,0.00,0.00}{342}1$}] {}
        -- ++(300:1.4cm) node (4321) [label=right: {\scriptsize$\textcolor[rgb]{0.00,0.00,0.00}{432}1$}] {}
        -- ++(240:1.4cm) node (2341) [label=135: {\scriptsize$\textcolor[rgb]{0.00,0.00,0.00}{234}1$}] {}
        -- ++(180:1.4cm) node (3241) [label=-70: {\scriptsize$\textcolor[rgb]{0.00,0.00,0.00}{324}1$}] {}
        -- ++(120:1.4cm) node (4231) [label=right:
        {\scriptsize$\textcolor[rgb]{0.00,0.00,0.00}{423}1$}] {}
        -- (2431); 

     \draw [line width = 1pt,red](-15,3.6) 
         ++(0:1.4cm) 
         ++(0:4.0cm) node (1423) [label=135: {\scriptsize$\textcolor[rgb]{0.00,0.00,0.00}{142}3$}] {}
       -- ++(0:1.4cm) node (4123) [label=45: {\scriptsize$\textcolor[rgb]{0.00,0.00,0.00}{412}3$}] {}
        -- ++(300:1.4cm) node (2143) [label=left: {\scriptsize$\textcolor[rgb]{0.00,0.00,0.00}{214}3$}] {}
        -- ++(240:1.4cm) node (1243) [label=-110: {\scriptsize$\textcolor[rgb]{0.00,0.00,0.00}{124}3$}] {}
        -- ++(180:1.4cm) node (4213) [label=45: {\scriptsize$\textcolor[rgb]{0.00,0.00,0.00}{421}3$}] {}
        -- ++(120:1.4cm) node (2413) [label=180:
        {\scriptsize$\textcolor[rgb]{0.00,0.00,0.00}{241}3$}] {}
        -- (1423); 
\draw [line width = 1pt,red](-15,3.6) 
          ++(240:1.4cm) 
         ++(270:5.0cm) node (1234) [label=right: {\scriptsize$\textcolor[rgb]{0.00,0.00,0.00}{123}4$}] {}
        -- ++(60:1.4cm) node (2134) [label=50: {\scriptsize$\textcolor[rgb]{0.00,0.00,0.00}{213}4$}] {}
       -- ++(0:1.4cm) node (3124) [label=225: {\scriptsize$\textcolor[rgb]{0.00,0.00,0.00}{312}4$}] {}
        -- ++(300:1.4cm) node (1324) [label=0: {\scriptsize$\textcolor[rgb]{0.00,0.00,0.00}{132}4$}] {}
        -- ++(240:1.4cm) node (2314) [label=-45: {\scriptsize$\textcolor[rgb]{0.00,0.00,0.00}{231}4$}] {}
        -- ++(180:1.4cm) node (3214) [label=-135:
        {\scriptsize$\textcolor[rgb]{0.00,0.00,0.00}{321}4$}] {}
       -- (1234); 
 \draw [line width = 1pt,red](-15,3.6) 
         ++(0:1.4cm) 
         ++(0:4.0cm) 
         ++(240:1.4cm) 
         ++(270:5.0cm) node (3412) [label=180: {\scriptsize$\textcolor[rgb]{0.00,0.00,0.00}{341}2$}] {}
          -- ++(60:1.4cm) node (1432) [label=300: {\scriptsize$\textcolor[rgb]{0.00,0.00,0.00}{143}2$}] {}
       -- ++(0:1.4cm) node (4132) [label=130: {\scriptsize$\textcolor[rgb]{0.00,0.00,0.00}{413}2$}] {}
       -- ++ (300:1.4cm)  node (3142) [label=180: {\scriptsize$\textcolor[rgb]{0.00,0.00,0.00}{314}2$}] {}
        -- ++(240:1.4cm) node (1342) [label=-45: {\scriptsize$\textcolor[rgb]{0.00,0.00,0.00}{134}2$}] {}
        -- ++(180:1.4cm) node (4312) [label=-135:
        {\scriptsize$\textcolor[rgb]{0.00,0.00,0.00}{431}2$}  ] {}
        -- (3412) ; 

    \draw [line width = 1pt,blue] (2431) -- (1432);
    \draw [line width = 1pt,blue](2341) -- (1342);
    \draw [line width = 1pt,blue] (4123) -- (3124);
    \draw [line width = 1pt,blue](4213) -- (3214);
    \draw [line width = 1pt,blue](4321) -- (1324);
    \draw [line width = 1pt,blue](2143) -- (3142);
    \draw [line width = 1pt,blue](2314) -- (4312);
    \draw [line width = 1pt,blue](2413) -- (3412);
    \draw [line width = 1pt,blue](3421) -- (1423);
    \draw [line width = 1pt,blue](4231) -- (1234);
    \draw [line width = 1pt,blue](3241) to [out=-20,in=200] (1243);
    \draw [line width = 1pt,blue](2134) to [out=20,in=-200] (4132);

\node [line width = 1.0pt](1342) at (30:5.5)      [label=10: {\scriptsize$\red{1}\black{342}$}]{};
\node [line width = 1.0pt](1324) at (1*60+30:5.5) [label=180: {\scriptsize$\red{1}\black{324}$}]{};
\node [line width = 1.0pt](1234) at (2*60+30:5.5) [label=10: {\scriptsize$\red{1}\black{234}$}]{};
\node [line width = 1.0pt](1243) at (3*60+30:5.5) [label=-10: {\scriptsize$\red{1}\black{243}$}]{};
\node [line width = 1.0pt](1423) at (4*60+30:5.5) [label=180: {\scriptsize$\red{1}\black{423}$}]{};

\node [line width = 1.0pt,white] at (4*60+30:3)
[label=-90: {\scriptsize$ Partitioned \ along\ symbol \ 1$}]{};

\node [line width = 1.0pt](1432) at (5*60+30:5.5) [label=-10: {\scriptsize$\red{1}\black{432}$}]{};

     \draw [line width = 1.pt,red]
     (120:1.4cm) node (2134)
    [label=-10: {\scriptsize$\textcolor[rgb]{0.00,0.00,0.00}{2\textcolor[rgb]{1,0.00,0.00}{1}34}$}]{}
      -- ++(0:1.4cm) node (4132)  [label=0: {\scriptsize$\textcolor[rgb]{0.00,0.00,0.00}{4\textcolor[rgb]{1,0.00,0.00}{1}32}$}]{}
   -- ++(300:1.4cm) node (3142) [label=-180: {\scriptsize$\textcolor[rgb]{0.00,0.00,0.00}{3\textcolor[rgb]{1.00,0.00,0.00}{1}42}$}] {}
     -- ++(240:1.4cm) node (2143)  [label=0: {\scriptsize$\textcolor[rgb]{0.00,0.00,0.00}{2\textcolor[rgb]{1.00,0.00,0.00}{1}43}$}]{}
     -- ++(180:1.4cm) node (4123)  [label=-180: {\scriptsize$\textcolor[rgb]{0.00,0.00,0.00}{4\textcolor[rgb]{1.00,0.00,0.00}{1}23}$}]{}
       -- ++(120:1.4cm) node (3124)  [label=0:
     {\scriptsize$\textcolor[rgb]{0.00,0.00,0.00}{3\textcolor[rgb]{1.00,0.00,0.00}{1}24}$}]{}
     -- (2134); 

 \draw [line width = 1.pt,red](120:3cm) node (2314)
         [label=-180: {\scriptsize$\textcolor[rgb]{0.00,0.00,0.00}{23\textcolor[rgb]{1.00,0.00,0.00}{1}4}$}] {}
        -- ++(0:3cm) node (4312) [ label=-170: {\scriptsize$\textcolor[rgb]{0.00,0.00,0.00}{43\textcolor[rgb]{1.00,0.00,0.00}{1}2}$}] {}
        -- ++(300:3cm) node (3412) [label=0: {\scriptsize$\textcolor[rgb]{0.00,0.00,0.00}{34\textcolor[rgb]{1.00,0.00,0.00}{1}2}$}] {}
        -- ++(240:3cm) node (2413) [label=170: {\scriptsize$\textcolor[rgb]{0.00,0.00,0.00}{24\textcolor[rgb]{1.00,0.00,0.00}{1}3}$}]{}
        -- ++(180:3cm) node (4213) [label=-70: {\scriptsize$\textcolor[rgb]{0.00,0.00,0.00}{42\textcolor[rgb]{1.00,0.00,0.00}{1}3}$}] {}
        -- ++(120:3cm) node (3214) [label=0:
        {\scriptsize$\textcolor[rgb]{0.00,0.00,0.00}{32\textcolor[rgb]{1.00,0.00,0.00}{1}4}$}] {}
        -- (2314); 

 \draw [line width = 1.pt,red](120:4.7cm) node (2341)
         [label=10: {\scriptsize$\textcolor[rgb]{0.00,0.00,0.00}{234}1$}] {}
        -- ++(0:4.7cm) node (4321)    [label=0: {\scriptsize$\textcolor[rgb]{0.00,0.00,0.00}{432}1$}] {}
        -- ++(300:4.7cm) node (3421)  [label=0: {\scriptsize$\textcolor[rgb]{0.00,0.00,0.00}{342}1$}] {}
        -- ++(240:4.7cm) node (2431)  [label=-10: {\scriptsize$\textcolor[rgb]{0.00,0.00,0.00}{243}1$}]{}
        -- ++(180:4.7cm) node (4231)  [label=-10: {\scriptsize$\textcolor[rgb]{0.00,0.00,0.00}{423}1$}] {}
        -- ++(120:4.7cm) node (3241)  [label=180:
        {\scriptsize$\textcolor[rgb]{0.00,0.00,0.00}{324}1$}] {}
        -- (2341); 

 \draw [line width = 1pt,blue] (1342) -- (2341);
 \draw [line width = 1pt,blue] (1342) -- (4312);
 \draw [line width = 1pt,blue] (1342) -- (3142);

 \draw [line width = 1pt,blue] (1324) -- (2314);
 \draw [line width = 1pt,blue] (1324) -- (3124);
 \draw [line width = 1pt,blue] (1324) -- (4321);

 \draw [line width = 1pt,blue] (1234) -- (2134);
 \draw [line width = 1pt,blue] (1234) -- (4231);
 \draw [line width = 1pt,blue] (1234) -- (3214);

 \draw [line width = 1pt,blue] (1243) -- (3241);
 \draw [line width = 1pt,blue] (1243) -- (4213);
 \draw [line width = 1pt,blue] (1243) -- (2143);

 \draw [line width = 1pt,blue] (1423) -- (4123);
 \draw [line width = 1pt,blue] (1423) -- (2413);
 \draw [line width = 1pt,blue] (1423) -- (3421);

 \draw [line width = 1pt,blue] (1432) -- (2431);
 \draw [line width = 1pt,blue] (1432) -- (3412);
 \draw [line width = 1pt,blue] (1432) -- (4132);
 \end{tikzpicture}
\vskip -2cm
\caption{\label{f1}\footnotesize {Two structures of the 4-dimensional star graph $S_4$}}

\end{center}
\vskip-16pt

\end{figure}
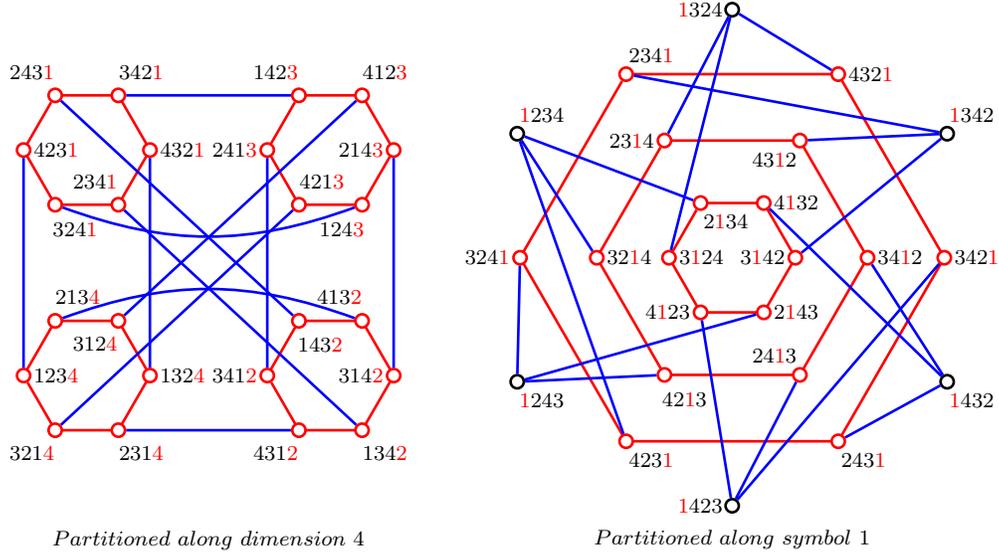


\vskip6pt

For a fixed symbol $i\in I_n$, let $S^{j:i}_n$ denote a subgraph of
$S_{n}$ induced by all vertices with symbol $i$ in the $j$-th
position for each $j\in I_n$. By the definition of $S_{n}$, it is
easy to see that $S^{j:i}_n$ is isomorphic to $S_{n-1}$ for each
$j\in I'_n$ and $S^{1:i}_n$ is an empty graph with $(n - 1)!$
vertices.

We will use two different hierarchical structures of $S_n$ depending
on different partition methods. The first one is partitioning along
a fixed dimension, which is clear and used frequently. The second
one is partitioning along a fixed symbol in $I_n$, which is a new
structure proposed recently by Shi~{\it et al.}\cite{sls12}.

\begin{lem}\label{lem2.1}{\rm (
The first structure)} For a fixed dimension $j\in I'_n$, $S_{n}$ can
be partitioned into $n$ subgraphs $S^{j:i}_n$, which is isomorphic
to $S_{n-1}$ for each $i\in I_n$. Moreover, there are $(n - 2)!$
independent edges between $S^{j:i_1}_n$ and $S^{j:i_2}_n$ for any
$i_1,i_2\in I_n$ with $i_1\ne i_2$.
\end{lem}

\begin{lem}\label{lem2.2}\textnormal{(Shi {\it et al.}~\cite{sls12}, 2012, The second structure)}
For a fixed symbol $i\in I_n$, $S_{n}$ can be partitioned into $n$
subgraphs $S^{j:i}_n$, which is isomorphic to $S_{n-1}$ for each
$j\in I'_n$ and $S^{1:i}_n$ is an empty graph with $(n - 1)!$
vertices. Moreover, there are a perfect matching  between
$S^{1:i}_n$ and $S^{j:i}_n$ for any $j\in I'_n$, and there are no
edge between $S^{j_1:i}_n$ and $S^{j_2:i}_n$ for any $j_1,j_2\in
I'_n$ with $j_1\ne j_2$.
\end{lem}

It is easy to know that $S_1,S_2,S_3$ are isomorphic to
$K_1,K_2,C_6$, respectively. $S_4$ is illustrated in Figure~\ref{f1}
by two different partition methods. As we will see, $S_4$ is
partitioned along dimension 4 in the left figure, and is partitioned
along symbol 1 in the right one.


\section{Main results}

In this section, we present our main results, that is, we determine
the $k$-super connectivity and $k$-super edge connectivity of the
$n$-dimensional star graph $S_{n}$. We first investigate the
properties of subgraphs in $S_n$ with minimum degree at least $k$.
For a subset $X\subseteq V(S_n)$ and $j\in I_n$, we use $U^X_j$ to
denote the set of symbols in the $j$-th position of vertices in $X$,
formally, $U^X_j=\{ p_j:\  p_{1}\ldots p_j\ldots p_{n}\in X \}$. The
following lemma plays a key role in our proof.

\begin{lem}\label{lem3.1}
Let $H$ be a subgraph of $S_n$ with vertex-set $X$ and $k\in
I_{n-1}$ a fixed integer. If $\delta(H)\geqslant k$, then there
exists some $j\in I'_n$ such that $|U^X_j|\geqslant k+1$.
\end{lem}

\begin{pf} Without loss of generality, we can assume that $H$ is connected. For sake of
simplicity, we write $U_j$ for $U^X_j$. Let $W_i$ be the set of
positions which symbol $i$ appears in vertices in $X$ excluding the
first position, that is, $W_i=\{j\in I'_n:\ i\in U_j\}$.

We use the second hierarchical structure of $S_n$ stated in
Lemma~\ref{lem2.2} to prove the lemma by induction on $n\,(\geqslant
k+1)$.

If $n=k+1$, then $\delta(H)\geqslant k=n-1$, and so $H=S_n$. Since
$|U_1|=|U_2|=\cdots=|U_n|=n=k+1$, the conclusion is hold for
$n=k+1$. We assume the conclusion is true for $n-1$ with $n\geqslant
k+2$.

Let $x=p_1p_2 \cdots p_n$ be a vertex in $H$. Then $x\in
V(S^{1:p_1}_n)$. By the second hierarchical structure, all the
neighbors of $x$ are in different $S^{j:p_1}_n$ for each $j\in
I'_n$.  Since $\delta(H)\geqslant k$, $p_1$ appears in at least $k$
different positions of vertices in $H$ excluding the first position.
It follows that
\begin{equation}\label{e3.1}
 |W_{p_1}|\geqslant k\ \ {\rm for\ any}\ x=p_1p_2 \cdots  p_n \in X
 \end{equation}

If $|U_1|=n$, then each symbol of $I_n$ appears in the first
position of vertices in $H$. By (\ref{e3.1}), we have
\begin{equation}\label{e3.2}
 |W_i|\geqslant k \ \ {\rm for\ each}\ i\in I_n.
 \end{equation}

Now we construct an $n\times (n-1)$ matrix $C=(c_{ij})_{n\times (n-1)}$,
where
 $$
 c_{ij}= \left\{\begin{array}{ll}
 1&{}\  j+1\in W_i\\ 
 0&{}\  {\rm otherwise}.\\
    \end{array}
\right.
 $$
Then
 $$
 \begin{array}{l}
 |U_j| =\sum\limits_{i=1}^n c_{ij}\ \ \text{for each $j\in I'_n$  and}\\
 |W_i| = \sum\limits_{j=2}^n c_{ij}\ \ \text{for each}\ i\in I_n.
 \end{array}
 $$
It follows that
\begin{equation}\label{e3.3}
\sum_{j=2}^n |U_j| =\sum_{j=2}^n \sum_{i=1}^n c_{ij}=\sum_{i=1}^n
\sum_{j=2}^n c_{ij}= \sum_{i=1}^n |W_i|.
  \end{equation}
Combining (\ref{e3.3}) with (\ref{e3.2}), we have
 \begin{equation}\label{e3.4}
 \sum_{j=2}^n |U_j|=\sum_{i=1}^n |W_i|\geqslant nk.
  \end{equation}
If $|U_j|\leqslant k$ for each $j\in I'_n$, then $(n-1)k\geqslant
nk$ by (\ref{e3.4}), a contradiction. Thus, there exists some $j\in
I'_n$ such that $|U_j|\geqslant k+1$.

If $|U_1|<n$, then there exists at least one symbol in $I_n$ that
does not appear in the first position of any vertex in $H$. Without
loss of generality, assume $1\not \in U_1$. Then $S^{1:1}_n$ does
not contain vertices of $H$. By the second hierarchical structure,
$H$ must be contained in the unique $S^{j_0:1}_n$ for some $j_0\in
I'_n$ since $H$ is connected. Because $S^{j_0:1}_n$ is isomorphic to
$S_{n-1}$, and $H\subseteq S^{j_0:1}_n$, by the induction
hypothesis, there exist some $j\in I'_n$ such that $|U_j|\geqslant
k+1$.

By the induction principle, the lemma follows.
\end{pf}

\begin{lem}\label{lem3.2}
For any integer $k$ with $0 \leqslant k \leqslant n-2$,
$\kappa_s^{(k)}(S_{n})\leqslant (k+1)!(n-k-1)$ and
$\lambda_s^{(k)}(S_{n})\leqslant (k+1)!(n-k-1)$.
\end{lem}

\begin{pf}
Let
 $$
 X=\{ \ p_1\cdots  p_{k+1}12\cdots (n-k-1):
 \  p_i \in I_n\setminus I_{n-k-1} { \rm \ for\ each}\ i \in  I_{k+1}\}.
 $$
Then, the subgraph $H$ of $S_n$ induced by $X$ is isomorphic to
$S_{k+1}$. Let $T$ be the set of neighbors of $X$ in $S_n-X$ and $F$
the set of edges between $X$ and $T$. By the definition of $S_n$,
 $$
 \begin{array}{rl}
 T=&\{i p_2 \cdots p_{k+1}12\cdots (i-1)p_1(i+1)\cdots (n-k-1):\ \\
  & \qquad \qquad i \in I_{n-k-1},\ p_j\in I_n\setminus I_{n-k-1}{ \rm \ for \ } j \ \in  I_{k+1}\}.
  \end{array}
 $$
For a vertex of $X$, since it has $k$ neighbors in $X$, it has
exactly $(n-k-1)$ neighbors in $T$. In addition, it is easy to see
that every vertex of $T$ has exactly one neighbor in $X$. It follows
that
 $$
  |T|=|F|=(k+1)!(n-k-1).
 $$

Since every vertex $v$ in $S_n-X$ has at most one neighbor in $X$,
$v$ has at least $((n-1)-1\geqslant)\, k$ neighbors in $S_n-X$,
which implies that $F$ is a $k$-edge-cut of $S_n$. It follows that
 $$
 \lambda_s^{(k)}(S_n)\leqslant |F|=(k+1)!(n-k-1).
 $$

We now show that $T$ is a $k$-vertex-cut of $S_n$. To this end, we
only need to show that every vertex in $S_n-(X\cup T)$ has at least
$k$ neighbors within.

Let $u$ be arbitrary vertex of $S_n-(X\cup T)$. We need to show that
at most one of neighbors of $u$ is in $T$. Suppose to the contrary
that $u$ has two distinct neighbors $v$ and $w$ in $T$. Then the
first digits of $v$ and $w$ are different. Without loss of
generality, assume $v=1 p_2 \ldots p_{k+1}p_123\cdots (n-k-1)$ and
$w=2 p'_2 \ldots p'_{k+1}1p'_23\cdots (n-k-1)$. Since $u$ is
adjacent to $v$, then $u$ and $v$ have exactly one digit difference
excluding the first one. So are $u$ and $w$. Therefore, $w$ and $v$
have exactly two digits difference excluding the first one. But $w$
and $v$ have yet two digits( the $(k+2)$-th and the $(k+3)$-th)
difference, then $p_2 \ldots p_{k+1}=p'_2 \ldots p'_{k+1}$,
therefore $v=w$, a contradiction.

Since $u$ has at most one neighbor in $T$, $u$ has at least
$((n-1)-1\geqslant) k$ neighbors in $S_n-(X\cup T)$, which implies
that $T$ is a $k$-vertex-cut of $S_n$. It follows that
 $$
 \kappa_s^{(k)}(S_n)\leqslant |T|=(k+1)!(n-k-1).
 $$
The lemma follows.
 \end{pf}

\vskip6pt

\begin{thm}\label{thm3.3}
$\kappa_s^{(k)}(S_{n})=\lambda_s^{(k)}(S_{n})=(k+1)!(n-k-1)$ for any
$k$ with $0 \leqslant k \leqslant n-2$.
\end{thm}

\begin{pf}
By Lemma~\ref{lem3.2}, we only need to show that, for any $k$ with
$0 \leqslant k \leqslant n-2$,
\begin{equation}\label{e3.5}
\lambda_s^{(k)}(S_{n})\geqslant(k+1)!(n-k-1) {\rm \ and \ }
\kappa_s^{(k)}(S_{n})\geqslant(k+1)!(n-k-1).
 \end{equation}

We prove (\ref{e3.5}) by induction on $k$. If $k=0$, then
$\lambda_s^{(0)}(S_{n})=\lambda(S_{n})=n-1$ and
$\kappa_s^{(0)}(S_{n})=\kappa(S_{n})=n-1$ by Lemma~\ref{lem2.3}, and
so (\ref{e3.5}) is true for $k=0$. Assume (\ref{e3.5}) holds for
$k-1$ with $k\geqslant 1$, that is, for any $k$ with $1\leqslant k
\leqslant n-2$,
\begin{equation}\label{e3.6}
\kappa_s^{(k-1)}(S_{n-1})\geqslant k!(n-k-1)\ \ {\rm and}\ \
\lambda_s^{(k-1)}(S_{n-1})\geqslant k!(n-k-1).
 \end{equation}

Let $T$ be a minimum $k$-vertex-cut (or $k$-edge-cut) of $S_{n}$. We
show that
 \begin{equation}\label{e3.7}
  |T|\geqslant (k+1)!(n-k-1)\ {\rm \ for}\ 1 \leqslant k \leqslant n-2.
 \end{equation}

 Let $X$ be the vertex-set of a connected
component $H$ of $S_{n}-T$, and
 $$
 Y=\left\{\begin{array}{ll}
 V(S_{n}-(X\cup T))\ & \text{if $T$ is a vertex-cut};\\
 V(S_{n}-X)\ & \text{if $T$ is an edge-cut}.
 \end{array}\right.
 $$
Then $\delta(H)\geqslant k$, and so there exists some $j\in I'_n$
such that $|U^X_{j}|\geqslant k+1$ by Lemma~\ref{lem3.1}. We choose
$j_0\in \{ j\in I'_n: \ |U^X_{j}|\geqslant k+1\}$ such that
$|U^X_{j_0}\cap U^Y_{j_0}|$ and $|U^Y_{j_0}|$ are as large as
possible. Without loss of generality, assume $j_0=n$. In the
following proof, we use the first hierarchical structure stated in
Lemma~\ref{lem2.1}. Let, for $i\in I_n$,
 $$
 \begin{array}{l}
 X_i=X\cap V(S^{n:i}_{n-1}), \ \ \  Y_i=Y\cap V(S^{n:i}_n),\\
 T_i=\left\{
 \begin{array}{rl}
 T\cap V(S^{n:i}_n) &\ \text{if  $T$ is a vertex-cut};\\
  T\cap E(S^{n:i}_n) &\ \text{if $T$ is an edge-cut},
 \end{array}\right.
  \end{array}
 $$
and let
 $$ J_X=\{i\in I_n:\ X_i\ne\emptyset\},\ \  J_Y=\{i\in I_n:\
Y_i\ne\emptyset\},
  \ J_0=J_X \cap J_Y.
  $$

Clearly, $ |J_X|=|U^X_{n}|, |J_Y|=|U^Y_{n}|$ and $|J_0|=|U^X_{n}\cap
U^Y_{n}|$.

If $i\in J_0$, $T_i$ is a vertex-cut (or an edge-cut) of
$S^{n:i}_n$. For any vertex $x$ in $S^{n:i}_n-T_i$, since $x$ has
degree at least $k$ in $S_{n}-T$ and has exactly one neighbor
outsider $S^{n:i}_n$, $x$ has degree at least $k-1$ in
$S^{n:i}_{n}-T_i$. Therefore, $T_i$ is a $(k-1)$-vertex-cut (or a
$(k-1)$-edge-cut) of $S^{n:i}_n$ for any $i\in J_0$. By the
induction hypothesis (\ref{e3.6}), we have
 \begin{equation}\label{e3.8}
 |T_i| \geqslant k!(n-k-1)\  {\rm \ for\ each}\ i \in J_0.
 \end{equation}

If $|J_0|\geqslant k+1$, by (\ref{e3.8}) we have
\[ |T|=\sum_{i=1}^n|T_i|\geqslant \sum_{i\in J_0}|T_i|\geqslant (k+1)k!(n-k-1)=(k+1)!(n-k-1),
\]
and so (\ref{e3.7}) follows.

Now assume $|J_0|\leqslant k$. Then $J_X\setminus J_0\ne
\emptyset$. We consider two cases, $J_Y\setminus J_0 \ne \emptyset $
and $J_Y\setminus J_0 =\emptyset$, respectively.

\vskip6pt

{\bf Case 1.}\  $J_Y\setminus J_0\ne  \emptyset $,

Assume $j_1\in J_X\setminus J_0,j_2\in J_Y\setminus J_0$. Then there
are $(n - 2)!$ independent edges between $S^{n:j_1}_n$ and
$S^{n:j_2}_n$. Since each vertex in $S^{n:j_1}_n$
 has unique external neighbor, thus $\bigcup_{j_1\in J_X\setminus J_0}S^{n:j_1}_n$
and $\bigcup_{j_2\in J_Y\setminus J_0}S^{n:j_2}_n$ have
$|J_X\setminus J_0||J_Y\setminus J_0|(n-2)!$ independent edges
between them. Note that each edge of which must have one end-vertex
in $T$ if $T$ is a vertex-cut, and each edge of which is contained
in $T$ if $T$ is an edge-cut. Therefore, no matter $T$ is a
vertex-cut or an edge-cut, we have
 \begin{equation}\label{e3.9}
 \sum_{i\in (J_X\cup J_Y)\setminus J_0}|T_i|\geqslant |J_X\setminus J_0||J_Y\setminus J_0|(n-2)!.
 \end{equation}
Let
 $$
 a=|J_X\setminus J_0|,\ b=|J_Y\setminus J_0|,\ c=|I_n\setminus (J_X\cup J_Y)|.
 $$
Then $a\geqslant 1, b\geqslant 1, a+b+c=n-|J_0|$, and so
 $$
 \begin{array}{rl}
 ab+c= & ab+(n-|J_0|)-(a+b)\\
     =& (n-|J_0|)+(a-1)(b-1)-1\\
     \geqslant & (n-|J_0|-1),
     \end{array}
 $$
that is,
 \begin{equation}\label{e3.10}
 ab+c\geqslant (n-|J_0|-1).
 \end{equation}

Note that $c=0$ if $T$ is an edge-cut. Thus if there exists some
$i\in I_n\setminus(J_X\cup J_Y)$, then $T$ is a vertex-cut and
$T_i=S^{n:i}_n$, and so
 \begin{equation}\label{e3.11}
 |T_i|=(n-1)! \ \ \text{if}\ i\in I_n\setminus(J_X\cup J_Y).
 \end{equation}
Thus, no matter $T$ is a vertex-cut or an edge-cut. Combining
(\ref{e3.8}), (\ref{e3.9}) and (\ref{e3.11}) with (\ref{e3.10}), we
have that
$$
\begin{array}{rl}
|T|&=\sum\limits_{i=1}^n|T_i| \geqslant \sum\limits_{i\in J_0}|T_i|
 + \sum\limits_{i\in (J_X\cup J_Y)\setminus J_0}|T_i|+\sum\limits_{i\in I_n\setminus(J_X\cup J_Y)}|T_i| \\
 &\geqslant |J_0|k!(n-k-1)+ ab(n-2)!+c(n-1)! \\
 &\geqslant |J_0|k!(n-k-1)+(ab+c)(n-2)! \\
 &\geqslant |J_0|k!(n-k-1)+(n-|J_0|-1)(n-2)! \\
 &\geqslant (n-1)k!(n-k-1)\\
 &\geqslant (k+1)!(n-k-1),
\end{array}
$$
and so (\ref{e3.7}) follows.

\vskip6pt

{\bf Case 2.}\ $J_Y\setminus J_0 = \emptyset $,

In this case $J_Y=J_0$, then $|U^Y_n|=|J_Y|\leqslant k$. Let
$\overline{X}_i=S^{n:i}_n-X_i$ for each $i\in I_n\setminus J_0$.
Note that for each $i\in I_n\setminus J_0$, $\overline{X}_i=T_i$ if
$T$ is a vertex-cut, and $\overline{X}_i=\emptyset$ if $T$ is an
edge-cut.

We first show there is no $i\in I_n \setminus J_0$ such that
$|\overline{X}_i|<(n-2)!$. Suppose to the contrary that there exists
some $i\in I_n\setminus J_0$ such that $|\overline{X}_i|<(n-2)!$.

We show $|U^{X_i}_j|\geqslant n-1$ for any $j\in I'_{n-1}$. On the
contrary, there exists some $j\in I'_{n-1}$ such that
$|U^{X_i}_j|\leqslant n-2$. Notice that the rightmost digit of every
vertex in $X_i$ is $i$. There is at least one symbol $i_1\in
I_n\setminus\{ i \}$ that does not appear in the $j$-th position of
any vertex in $X_i$. Thus, the vertices with symbol $i_1$ in the
$j$-th position and symbol $i$ in the $n$-th position are not
contained in $X_i$, which means that $\overline{X}_i$ contains at
least $(n-2)!$ vertices, that is, $|\overline{X}_i|\geqslant
(n-2)!$, a contradiction. Thus, $|U^{X_i}_j|\geqslant n-1$, and so
$|U^{X}_j|\geqslant n-1$ for any $j\in I'_{n-1}$.

Since $\delta(Y)\geqslant k$ and $|U_n^Y|\leqslant k$, by
Lemma~\ref{lem3.1} there exists some $j_1\in I'_{n-1}$ such that
$|U^Y_{j_1}|\geqslant k+1$. Then $|U^{X}_{j_1}|\geqslant n-1$ and
$|U^Y_{j_1}|\geqslant k+1$, and so $|U^{X}_{j_1}\cap
U^{Y}_{j_1}|\geqslant k$ and $|U^Y_{j_1}|\geqslant k+1$. Note that
$|U^{X}_{n}\cap U^{Y}_{n}|=|J_0|\leqslant k$ and $|U^Y_{n}|=|J_Y| =
|J_0|\leqslant k$. This contradicts to the choice of $j_0\, (=n)$.

Thus, there is no $i\in I_n\setminus J_0$ such that
$|\overline{X}_i|<(n-2)!$, and so there is no $i\in I_n\setminus
J_0$ such that $|\overline{X}_i|=0$. If $T$ is an edge-cut, then
$\overline{X}_i=\emptyset$, a contradiction. Therefore, $T$ is a
vertex-cut, and so $\overline{X}_i=T_i$. It follows that
 \begin{equation}\label{e3.12}
 |T_i|=|\overline{X}_i|\geqslant (n-2)!\ \ \text{for each}\ i\in I_n\setminus J_0.
 \end{equation}
Combining (\ref{e3.12}) with (\ref{e3.8}), we have
 $$
\begin{array}{rl}
|T|&=\sum\limits_{i=1}^n|T_i|=\sum\limits_{i\in J_0}|T_i|+\sum\limits_{i\in I_n\setminus J_0}|T_i|\\
&\geqslant |J_0|k!(n-k-1)+(n-|J_0|)(n-2)! \\
 &\geqslant (k+1)!(n-k-1).
\end{array}
$$
By induction principles, (\ref{e3.7}) holds and so the theorem
follows.
\end{pf}

\begin{cor}\textnormal{(~\cite{wz09}, ~\cite{ylm10})}
$\kappa_s^{(2)}(S_{n})=\lambda_s^{(2)}(S_n)=6(n-3)$ for $n\geqslant 4$.
\end{cor}

\section{Conclusions}

In this paper, we consider the generalized measures of fault
tolerance for networks, called the $k$-super connectivity
$\kappa_s^{(k)}$ and the $k$-super edge-connectivity
$\lambda_s^{(k)}$. For $n$-dimensional star graph $S_n$, which is an
attractive alternative network to hypercubes, we prove that
$\kappa_s^{(k)}(S_n)=\lambda_s^{(k)}(S_n)=(k+1)!(n-k-1)$ for $0 \leqslant
k \leqslant n-2$, which gives an affirmative answer to the
conjecture proposed by Wan and Zhang~\cite{wz09}. The results show
that at least $(k+1)!(n-k-1)$ vertices or edges have to be removed
from $S_n$ to make it disconnected and no vertices of degree less
than $k$. Thus these results can provide more accurate measurements
for fault tolerance of the system when $n$-dimensional star graphs
is used to model the topological structure of a large-scale parallel
processing system.


\vskip6pt


\begin{thebibliography}{10}

\bibitem{ak89}S. B. Akers, B. Krishnamurthy, A group theoretic model for symmetric
interconnection networks. IEEE Transactions on Computers, 38 (4)
(1989), 555-566.

\bibitem{cl02}
E. Cheng, M. J. Lipman, Increasing the connectivity of the star
graphs.  Networks, 40 (3) (2002), 165-169.

\bibitem{dt94}
K. Day, A. Triphthi, A comparative study of topological properties
of hypercubes and star graphs. IEEE Transactions on Parallel and
Distributed Systems, 5 (1) (1994), 31-38.


\bibitem{E89}
A. H. Esfahanian, Generalized measures of fault tolerance with
application to $n$-cube networks. IEEE Transactions on Computers, 38
(11) (1989), 1586-1591.

\bibitem{hy97}
S.-C. Hu, C.-B. Yang, Fault tolerance on star graphs. International
Journal of Foundations of Computer Science,  8 (2)(1997), 127-142.

\bibitem{lhm94}
S. Latifi, M. Hegde, M. Naraghi-Pour, Conditional connectivity
measures for large multiprocessor systems. IEEE Transactions on
Computers, 43 (2) (1994), 218-222.

\bibitem{sls12}%
W. Shi, F. Luo, P. K. Srimani, A new hierarchical structure of star
graphs and applications. In: Distributed Computing and Internet
Technology, LNCS, 7154, (2012), pp. 267-268.

\bibitem{rls96}%
Y. Rouskov, S. Latifi, P. K. Srimani, Conditional fault diameter
of star graph networks. Journal of Parallel and Distributed
Computing, 33 (1) (1996), 91-97.

\bibitem{wl10}
D. Walker, S. Latifi, Improving bounds on link failure tolerance of
the star graph. Information Sciences, 180 (13) (2010), 2571-2575.

\bibitem{wz09}
M. Wan, Z. Zhang, A kind of conditional vertex connectivity of star
graphs. Applied Mathematics Letters, 22 (2009), 264-267.

\bibitem{ylm10}
W.-H. Yang, H.-Z Li, J.-X Meng, Conditional connectivity of Cayley
graphs generated by transposition trees. Information Processing
Letters, 110 (23) (2010), 1027-1030.

\end{thebibliography}
\end{document}